\documentclass[runningheads,12pt,a4]{article}
\usepackage[left=2.5cm,right=2.5cm,top=3cm,bottom=3cm]{geometry}
\usepackage{amsmath,amssymb,amsfonts,amsthm}
\newcommand{\ecc}{\mathrm{ecc}}
\newcommand{\diam}{\mathrm{diam}}
\newcommand{\rad}{\mathrm{rad}}
\newtheorem{theorem}{Theorem}
\newtheorem{proposition}{Proposition}
\newtheorem{lemma}{Lemma}
\newtheorem{example}{Example}
\newtheorem{definition}{Definition}

\usepackage{setspace}
\usepackage{float}
\usepackage{authblk}
\date{}
\usepackage{tikz}
\begin{document}
	\title{Boundary vertices of Strongly Connected Digraphs with respect to `Sum Metric'}
 \author[1]{Bijo S. Anand}
 \author[2]{Manoj Changat}
\author[3]{Prasanth G. Narasimha-Shenoi}
\author[3]{Mary Shalet Thottungal Joseph}
\author[4]{Mithra R}
\author[5]{Prakash G. Narasimha-Shenoi}
\affil[1]{Department of Mathematics, Sree Narayana College, Punalur, Kollam, Kerala, India-691306\\bijos$\_$anand@yahoo.com}
\affil[2]{Department of Futures Studies, University of
    Kerala,Thiruvananthapuram, Kerala - 695 581, India\\mchangat@gmail.com}
\affil[3,4]{Department of Mathematics\\ Government College Chittur\\ Palakkad, India - 678104\\prasanthgns@gmail.com, mary$\_$shallet@yahoo.co.in, mithrar1729@gmail.com}
\affil[5]{Department of Mathematics\\ Maharajas College Ernakulam\\ Kerala, India - 682011\\prakashgn@gmail.com}
\affil[3,5]{Department of Collegiate Education, Government of Kerala, \\Thiruvananthapuram, Kerala, India - 695033}




	\maketitle 	
	\begin{abstract}
		Suppose $D = (V, E)$ is a strongly connected digraph and $u, v \in V (D)$. 
        Among the many metrics in graphs, the sum metric warrants further exploration. The sum
distance $sd(u, v)$ defined as $sd(u, v) =\overrightarrow{d}(u, v)+\overrightarrow{d}(v, u)$ is a metric where $\overrightarrow{d}(u, v)$
denotes the length of the shortest directed $u - v$ path in $D$. The four main boundary vertices
in the digraphs are “boundary vertices, contour vertices, eccentric vertices”, and
“peripheral vertices” and their relationships have been studied. Also, an attempt is
made to study the boundary-type sets of corona product of (di)graphs.  The center of the corona product of two strongly connected digraphs is established.  All the boundary-type sets and the center of the corona product are established in terms of factor digraphs.\\
  \noindent
	\textbf{Keywords:} Sum distance. Boundary-type sets. Strongly connected digraph. Corona product. \\
	\noindent	
		\textbf{Mathematics Subject Classification:} 05C12, 05C20, 05C76.
	\end{abstract}
	\section{Introduction}
	Graph theory is one of the popular areas in discrete mathematics with numerous theoretical developments and countless applications in practical problems. Its origin dates back to 1736 with the solution of the famous K\"{o}nigsberg Bridge Problem by the mathematician Leonhard Euler. Graphs can be partitioned into two types such as undirected graphs and directed graphs (digraphs). Even though both
		have enormous important applications, on various grounds, undirected	graphs have been studied much more widely than directed graphs. The undirected graphs form a special class of directed graphs (symmetric digraphs). Digraphs have a remarkable role in all aspects of scientific study. Some of the initial works in the discipline of digraphs are \cite{barton1968structural,bender1988asymptotic,roberts1975signed}.\\		
		The distance between two vertices $u$ and $v$ in a connected graph $G$ is the length of the shortest $u-v$ path in $G$. The directed distance of a digraph $D$ is the length of the shortest path from $u$ to $v$, where $u$ and $v$ are vertices in $D$, but it can be seen that this is not a metric.  Chartrand and Tian introduced two distances in strong digraphs, namely maximum distance and sum distance in the article `Distance in digraphs' \cite{chartrand1997distance}. Chartrand et al. have also defined a metric strong distance in the article `Strong distance in strong digraphs' \cite{chartrand1999strong}. Bielak and Syslo studied the peripheral
vertices in \cite{bielak1983peripheral}, and eccentric vertices by Chartrand et al. in \cite{chartrand1996eccentric}. Chartrand et al. defined a new boundary-type set called the boundary of a graph and also studied the relationship between the periphery, eccentricity, and boundary \cite{chartrand2003boundary}. The contour set was defined by C\'aceres et al. in \cite{caceres2005rebuilding}.  The center and radius of a graph have been a subject of study since \cite{kopylov1977centers}. The boundary-type sets of a graph namely boundary, contour, eccentricity and, periphery sets have also been studied in \cite{caceres2005rebuilding}.  The sum distance spectral properties of digraphs have been studied by \cite{xi2023general,xu2022extremal}
\\
 For the sake of completion, the boundary-type sets of the corona product of simple graphs and similar studies on the corona product of two directed graphs have also been discussed.
	\section{Preliminaries}
	A \textit{graph} $G=(V, E)$ consists of a nonempty set $V(G)=V$ called the vertex set and another set $E(G)=E$ whose elements are called edges such that each edge $e_k$ is identified with an unordered pair $(v_i,v_j)$ of vertices \cite{bondy1976graph}. An edge having the same vertex as both its end vertices is called a \textit{loop}. A set of two or more edges of a graph $G$ is called \textit{parallel} edges if they have the same distinct end vertices. A \textit{simple graph} is a graph with no loops and parallel edges. 
	
	A directed graph or \textit{digraph} $D$ consists of a non-empty finite set $V(D)$ of elements called vertices and a finite set $E(D)$ of ordered pairs of distinct vertices called edges \cite{bang2008digraphs}.  A digraph $D$ is denoted as $D=(V, E)$ where $V(D)=V$ is the vertex set and $E(D)=E$ is the edge set. For an edge $(u, v)$ the first vertex $u$ is its tail and the second vertex $v$ is its head. The definition of a digraph allows to have edges with the same end-vertices but does not allow it to contain parallel edges, that is, pairs of edges with the same tail and the same head, or loops (edges whose head and tail coincide) \cite{bang2008digraphs}. Let $D=(V, E)$ be a digraph and $u\in V$. The \textit{out-neighbourhood} of $u$, $N^+(u)=\{v\in V(D): (u, v) \in E(D)\}$ and the \textit{in-neighbourhood} of $u$, $N^-(u)=\{v\in V(D): (v,u) \in E(D)\}$. The out-neighbours and in-neighbours together form the neighbours of $u$ and this set is denoted as $N(u)$. Thus $N(u)=N^+(u) \bigcup N^-(u)$. The directed path in a directed graph can be described as a sequence of distinct vertices and a directed edge where the edge is pointing from each vertex in the sequence to its successor in the sequence.	A digraph $D$ is strongly connected or strong if for every pair $u, v$ of distinct vertices in $D$, there exists a path from $u$ to $v$. A digraph is said to be \textit{weakly connected} if there is a path between every pair of vertices when the direction of the edges is ignored.
	\subsection{Distance in Digraphs}
	The (directed) distance $\overrightarrow{d}(u, v)$ from $u$ to $v$ in a strong digraph $D$ is the length of the shortest $u-v$ path in $D$ where the length of a path in $D$ is the number of edges in the path  \cite{chartrand1997distance}. A $u-v$ path of length $\overrightarrow{d}(u,v)$ is called a $u-v$ geodesic \cite{bang2008digraphs}. But this distance is not a metric since $\overrightarrow{d}(u, v) \neq \overrightarrow{d}(v, u)$ in general. Hence two metrics on strong digraphs were introduced by
	Chartrand and Tian \cite{chartrand1997distance} namely \textit{maximum distance} $md(u, v)=\max\{\overrightarrow{d}(u, v), \overrightarrow{d}(v, u)\}$ and \textit{sum distance} $sd(u, v)=\overrightarrow{d}(u, v)+\overrightarrow{d}(v, u)$. Here we deal with the second metric, the sum distance $sd$. Even though Chartrand et al. proved in \cite{chartrand1997distance} that this sum distance is a metric, for the sake of completeness we prove that ``the Sum distance" is a metric, see as stated in Proposition \ref{sdmetric}.
	\begin{proposition}\label{sdmetric}
		If $D$ is a strong digraph, then the sum distance (sd) defined on $D$ is a metric.
	\end{proposition}
	\begin{proof}
		Let $D$ be a strong digraph. Suppose $u, v \in V(D)$. Then $sd(u, v)=\overrightarrow{d}(u, v)+\overrightarrow{d}(v, u)>0$ for all vertices $u$ and $v$. Suppose $sd(u, v)=0$. Then 
		\begin{align*}
			sd(u, v)=0 & \iff \overrightarrow{d}(u, v)+\overrightarrow{d}(v, u)=0\\
			& \iff \overrightarrow{d}(u, v)=0 \text{ and }  \overrightarrow{d}(v, u)=0\\
			& \iff u=v.\\
			sd(u,v)&= \overrightarrow{d}(u, v)+\overrightarrow{d}(v, u)\\
			&=\overrightarrow{d}(v, u)+\overrightarrow{d}(u, v)\\
			&=sd(v, u)\\
			sd(u, v)&=sd(v, u).
		\end{align*} 
		Hence $sd$ is positive and symmetric. Now to show that $sd$ satisfy the triangle inequality. Let $u, v, w \in V(D)$. 
		\begin{align*}
			sd(u, v)&=\overrightarrow{d}(u, v)+\overrightarrow{d}(v, u)\\
			&\leq \overrightarrow{d}(u, w)+\overrightarrow{d}(w, v)+\overrightarrow{d}(v, w)+\overrightarrow{d}(w, u)\\
			&= \overrightarrow{d}(u, w)+\overrightarrow{d}(w, u)+\overrightarrow{d}(w, v)+\overrightarrow{d}(v, w)\\
			&= sd(u, w)+sd(w, v)\\
			sd(u, v)&\le sd(u, w)+sd(w, v).
		\end{align*}
		Hence $sd$ is a metric.   
	\end{proof}
	The geodesic interval between two vertices $u$ and $v$ in $D$ is denoted by $I(u, v)$, the set of vertices of all shortest paths between $u$ and $v$ or the set of vertices of all shortest paths between $v$ and $u$ in $D$ \cite{caceres2006geodetic}. Let $S$ be a set of vertices of $D$. The geodetic closure $I[S]$ is the union of all geodetic intervals $I[u, v]$ between all pairs of vertices $u, v \in S$ \cite{caceres2006geodetic}. Thus $I[S]=\bigcup\limits_{u, v\in S} I[u, v]$. A set $S$ of vertices of $D$ is a geodetic set in $D$ if $I[S]=V(D)$. 
	\begin{example}
	Consider the first digraph in Figure \ref{fig2}. The vertices in the shortest path from $u_1$ to $u_4$ are $u_1,u_6, u_4$ and the vertices in the shortest path from $u_4$ to $u_1$ are $u_4,u_5,u_6,u_2,u_1$. Then $I[u_1,u_4]= \{u_1, u_2, u_4, u_5, u_6\}$.
 	\begin{figure}[H]
		\centering
		\resizebox{10cm}{6cm}{%
			\begin{tikzpicture}
				\tikzstyle{vertex} = [circle,fill=white,draw,minimum size=.65cm,inner sep=2pt,thick]
				\tikzstyle{label} = [circle,fill=white]
				\tikzstyle{label} = [fill=white]
				\tikzstyle{ecc} = [fill=white,thick]
				\tikzstyle{edge} = [->, thick]
				\tikzstyle{2edge} = [<->, thick]
				\node[vertex] (1) at (-2,2){$u_{1}$};
				\node[ecc](a) at(-2.7,2){$\textcolor{red}{\textbf{6}}$};
				\node[vertex] (2) at (0,4){$u_{2}$};		
				\node[ecc](b) at(0,4.7){$\textcolor{red}{\textbf{5}}$};	
				\node[vertex] (3) at (2,2){$u_{3}$};
				\node[ecc](c) at(2.7,2){$\textcolor{red}{\textbf{4}}$};
				\node[vertex] (4) at (2,-2){$u_{4}$};
				\node[ecc](d) at(2.7,-2){$\textcolor{red}{\textbf{6}}$};
				\node[vertex] (5) at (-2,-2){$u_{5}$};
				\node[ecc](e) at(-2.7,-2){$\textcolor{red}{\textbf{4}}$};
				\node[vertex] (6) at (0,0){$u_{6}$};
				\node[ecc] (f) at(0,-0.7){$\textcolor{red}{\textbf{4}}$};
				\node[vertex](7) at (7,-2) {$v_{1}$};
				\node[ecc](g) at(6.4,-2){$\textcolor{red}{\textbf{7}}$};
				\node[vertex](8) at (6,1) {$v_{2}$};
				\node[ecc](h) at(5.4,1){$\textcolor{red}{\textbf{5}}$};
				\node[vertex](9) at (6,4) {$v_{3}$};
				\node[ecc](i) at(5.4,4){$\textcolor{red}{\textbf{4}}$};
				\node[vertex](10) at (8,4) {$v_{4}$};
				\node[ecc](j) at(8.6,4){$\textcolor{red}{\textbf{7}}$};
				\node[vertex](11) at (8,1) {$v_{5}$};
				\node[ecc](k) at(8.6,1){$\textcolor{red}{\textbf{4}}$};
				\node[label] (12) at (0,-3){(1)};
				\node[label] (13)at (7,-3){(2)};
				\draw[2edge](1)--(2);
				\draw[2edge](2)--(3);
				\draw[edge](3)--(4);
				\draw[edge](4)--(5);
				\draw[edge](5)--(6);
				\draw[edge](1)--(5);
				\draw[edge](1)--(6);
				\draw[edge](6)--(2);
				\draw[edge](6)--(3);
				\draw[edge](6)--(4);
				\draw[2edge](7)--(8);
				\draw[2edge](8)--(9);
				\draw[edge](9)--(10);
				\draw[edge](10)--(11);
				\draw[edge](7)--(11);
				\draw[edge](11)--(9);
			\end{tikzpicture}		
		}
		\caption{}\label{fig2}
	\end{figure}
	\end{example}
	The Boundary, Contour, Eccentricity and Periphery sets of a strong digraph are discussed with respect to the metric. Hereafter $sd(u, v)$ is denoted as $d(u, v)$.
	
	\begin{definition} \cite{chartrand1997distance}
		The radius $rad(D)$ is the minimum eccentricity of vertices.
	\end{definition}
	\begin{definition}\cite{chartrand1997distance}
	The diameter $diam(D)$ is the maximum eccentricity of vertices.
	\end{definition}
	\begin{definition}\cite{chartrand1997distance}
		A vertex $u\in V(D)$ is called a central vertex of $D$ if its eccentricity is minimum over all vertices in $D$, that is, if the eccentricity of $u$ is equal to the radius of $D$. The center $C(D)$ of a strong digraph $D$ is the set of all central vertices.  That is
		\begin{align*}
			 C(D)&=\{v \in V(D): \ecc(v) \le \ecc(u), \forall u \in V(D)\}\\
			&=\{v \in V(D):\ecc(v)=\rad(D)\}.
		\end{align*}
	\end{definition}
	\subsection{Boundary-type Sets}
	The following definitions are analogous to the definitions in Chartrand et al.\cite{chartrand2003boundary}.
	\begin{definition}\cite{chartrand2003boundary}.
		Let $D = (V, E)$ be a strong digraph and $u, v \in V(D)$. The vertex $v$ is said to be a boundary vertex of $u$ if no neighbor of $v$ is further away from $u$ than $v$. A vertex $v$ is called a boundary vertex of $D$ if it is the boundary vertex of some vertex $u \in V(D)$. The boundary $\partial(D)$ of $D$ is the set of all its boundary vertices.  That is, $\partial(D)=\{v\in V(D): \exists\, u \in V(D), \forall w \in N(v):d(u, w)\leq d(u, v)\}$.
	\end{definition}
	\begin{definition}\cite{chartrand2003boundary}.	The eccentricity of a vertex $u \in D$ is defined as $\ecc_{D}(u)=\max\{d(u, v)|v\in V(D)\}$. If $D$ is clear from the context, then $\ecc_{D}(u)$ is denoted as $\ecc(u)$. Let $u, v \in V(D)$, the vertex $v$ is called an eccentric vertex of $u$ if no vertex in $V(D)$ is further away from $u$ than $v$, that is, if
		$d(u, v) = \ecc(u)$. A vertex $v$ is called an eccentric vertex of $D$ if it is the eccentric vertex of some vertex $u \in V(D)$. The eccentricity $Ecc(D)$ of $D$ is the set of all of its eccentric vertices. i.e; $Ecc(D) =\{v\in V(D):\exists u\in V(D), \ecc(u)=d(u, v)\}$.
	\end{definition}
	\begin{definition}\cite{chartrand2003boundary}. A vertex $v \in V(D)$ is called a peripheral vertex of $D$ if no vertex in $V(D)$ has an eccentricity greater than $ecc(v)$, that is, if the eccentricity of $v$ is exactly equal to the diameter $diam(D)$ of $D$. The periphery $Per(D)$ of $D$ is the set of all of its peripheral vertices.
		\begin{align*}
			Per(D)&=\{v\in V(D):\ecc(u)\le \ecc(v), \forall u \in V(D)\}\\
			&=\{v\in V(D): \ecc(v)=\diam(D)\}
		\end{align*}
	\end{definition}
	\begin{definition}\cite{chartrand2003boundary}. A vertex $v\in V(D)$ is called a contour vertex of $D$ if no neighbor vertex of $v$ has an eccentricity greater than $ecc(v)$. The contour $Ct(D)$ of $D$ is the set of all of its contour vertices.  That is $Ct(D)=\{v\in V(D): \ecc(u)\le \ecc(v), \forall u\in N(v)\}$.
	\end{definition}	
	Similar to the Proposition 6 in \cite{caceres2006geodetic} by C\'aceres et al in the case of simple graphs, we here prove the case of strong digraphs when the sum distance is considered.
	\begin{proposition}
		Let $D=(V, E)$ be a strong digraph. Then the following statements hold.	
		\begin{enumerate}
			\item $Per(D) \subseteq Ct(D) \bigcap Ecc(D)$.
			\item $Ecc(D) \bigcup Ct(D) \subseteq \partial(D)$.
		\end{enumerate} 
	\end{proposition}
	\begin{proof}
		Given $D=(V,E)$ be a strong digraph.
  \begin{enumerate}
      \item  Let $v \in Per(D)$. Therefore, for every $u \in V(D)$ we  have$\ecc(u) \le \ecc(v)$. So $\ecc(u) \le \ecc(v)$ for every $u \in N(v)$ and $v \in Ct(D)$. Now it remains to prove that $Per(D) \subseteq Ecc(D)$. Let $ w \in V(D) $ such that $\ecc(v)=d(v, w)$. Then $w \in Ecc(D)$. Since $w\in V(D)$ and $v\in Per(D)$, $\ecc(w)\le \ecc(v)$. Thus $\ecc(w) \le d(v,w)$. Also $\ecc(w)=\max\{d(w, u)|u \in V(D)\}$. Thus $d(w, v) \le \ecc(w)$. Hence $\ecc(w) = d(w,v)$. Therefore $v \in Ecc(D)$. Hence $Per(D) \subseteq Ct(D) \bigcap Ecc(D)$.
	\item  Let $v \in Ecc(D)\bigcup Ct(D)$. Then either $v \in Ecc(D)$ or $v \in Ct(D)$.  Suppose $v \in Ecc(D)$.	Then there exists $u \in V(D)$ such that $\ecc(u)=d(u,v)$. Now let $t \in N(v)$. Since $\ecc(u) = \max\{d(u,w):w \in V(D)\}$, $d(u,v) \geq d(u,w)$ for every $w\in V(D)$. This implies $d(u,v) \geq d(u,t) \text{ for every } t\in N(v)$. Thus $v$ is a boundary vertex of $u$. Therefore $v \in   \partial(D)$. Suppose $v \in Ct(D)$. Then $\ecc(u)\leq ecc(v)$,  for every $u \in N(v)$. Let $\ecc(v)=d(v,w)$. Then $\ecc(u) \leq d(v,w)$. Also $d(u,w) \leq d(v,w), \text{ for every } u \in N(v)$.
		Thus $v$ is a boundary vertex of $w$. Hence $v \in \partial(D)$. Thus  $Ecc(D) \bigcup Ct(D) \subseteq \partial(D)$.
  \end{enumerate}
	\end{proof}

	\section{Corona Product of Undirected Simple Graphs}
		Let $G$ be a connected graph of order $n$ and $H$ (connected) be any graph with $|V(H)|\geq 2$. Then the corona product $G\odot H$ of $G$ and $H$ is defined as a graph which is formed by taking $n$ copies $H_1,H_2, \ldots, H_n$ of $H$ and connecting  $i^{th}$
		vertex of $G$ to the vertices of $H_i$. Throughout this section $H_i$ is referred to 
		as  $i^{th}$ copy of $H$ connected to  $i^{th}$ vertex of $G$ in $G\odot H$ for every $i\in \{1,2, \ldots ,n\}$ \cite{iswadi2011metric}.
	Let $V(G)=\{u_1,u_2, \ldots,u_n\}$ and $V(H)=\{v_1,v_2, \ldots,v_m\}$ and $H_i, i\in I=\{1, 2,  \ldots, n\}$ be the $i^{th}$ copy of $H$. Each vertex $v_r$ of a copy of $H$ attached to a vertex $u_i$ can be uniquely described as $v_r^i$. Thus $V(H_i)=\{v_1^i,v_2^i, \ldots,v_m^i\}$. Corona product of $G$ and $H$ is a graph $G\odot H$ such that $V(G\odot H)=V(G)\bigcup \bigg\{\bigcup\limits_{i=1}^nV(H_i)\bigg\}$ and $E(G\odot H)=E(G)\bigcup \bigg\{\bigcup\limits_{i=1}^n E(H_i)\bigg\} \bigcup \bigg\{\bigcup \limits_{i=1}^n(u_i,v_r^i), r=1,2, \ldots,m\bigg\}$.\\
	
	\subsection{Distance Between Two Vertices in Corona Product of Graphs}
	The following proposition gives the distance between two vertices in Corona product of Simple Graphs. For the sake of completion,  we will prove the proposition.
	\begin{proposition} \label{distcorgraph}
		Let $G$ and $H$ be two simple graphs with vertex sets $\{u_1,u_2,  \ldots ,u_n\}$ and $\{v_1,v_2,  \ldots ,v_m\}$ respectively. Then
		\begin{equation*}
			d_{G \odot H}(x,y)  = \begin{cases}
				d_G(u_i,u_j) &\quad \text{if } x=u_i, y= u_j\\
				d_G(u_i,u_j)+1  &\quad \text{if } x=u_i, y=v_t^j \\
				d_G(u_i,u_j)+2 &\quad \text{if }x=v_r^i, y=v_t^j,i\neq j \\
				\min\{d_H(v_r,v_t),2\} &\quad \text{if }x=v_r^i, y=v_t^j,i=j.
			\end{cases}	
		\end{equation*}
	\end{proposition}
	\begin{proof}
		Let $G$ and $H$ be two graphs. Consider the vertices $ x,y \in V(G \odot H)$.\\
		\textbf{Case 1:} $x=u_i, y= u_j$ where $x \text{ and }y $ are vertices of $G$.\\
		Then $d_{G\odot H}(x,y)=d_{G\odot H}(u_i, u_j)$. Thus $d_{G\odot H}(x,y)= d_G(u_i,u_j) \text{ since } u_i \text{ and } u_j \text{ are vertices of } G$.
		Therefore $d_{G\odot H}(u_i,u_j)=d_G(u_i,u_j)$.\\
		\textbf{Case 2:} $x=u_i, y=v_t^j$.
		\begin{align*}
			d_{G\odot H}(x,y)&=d_{G\odot H}(u_i,v_t^j)\\
			&=d_{G\odot H}(u_i,u_j)+d_{G\odot H}(u_j,v_t^j)\\
			&=d_G(u_i,u_j)+1.
		\end{align*}
		\textbf{Case 3:} $x=v_r^i,y=v_t^j,i\neq j$.\\
		Now	$d_{G \odot H}(v_r^i,v_t^j)=d_{G\odot H}(v_r^i,u_i)+d_{G \odot H}(u_i,u_j)+d_{G\odot H}(u_j,v_t^j)$. Then $d_{G \odot H}(v_r^i,v_t^j)=d_G(u_i,u_j)+2$.\\
		\textbf{Case 4:} $x=v_r^i, y=v_t^j,i=j $.\\
  In this case the vertices $v_r^i$ and  $v_t^j$ are in the same component $H^i$.  If they are adjacent then $d_{G \odot H}(v_r^i,v_t^j)=1$  and if they are not adjacent  then the path $v_r^iu_iv_t^j$ is a path of length $2$ and so in this case the length of the shortest $v_r^i-v_t^j$ path is $2$.  Combining these two we have $d_{G\odot H}(v_r^i,v_t^i) = \min\{d_H(v_r,v_t),2\} \text{ when } i=j$.  Hence the proposition.
	\end{proof}
	
	\begin{lemma} \label{eccradius} \cite{yarahmadi2014center}
		Let $G$ and $H$ be graphs. Then
		\begin{enumerate}
			\item 	$\ecc_{G\odot H}(x)= \begin{cases}
				\ecc_G(u_i)+1, \quad x= u_i\\
				\ecc_G(u_i)+2,  \quad x= v_r^i.		
			\end{cases}$
			\item $r(G\odot H)=r(G)+1$.
			\item $d(G\odot H)=d(G)+2$.
		\end{enumerate}
	\end{lemma}
	
	\begin{theorem}\cite{yarahmadi2014center}
		Let $G$ and $H$ be two graphs. Then
		\begin{enumerate}
			\item $C(G\odot H) = C(G)$.
			\item $Per(G\odot H)=\{v_r^i:u_i\in Per(G)\}$.
		\end{enumerate}
	\end{theorem}
	\subsection{Boundary-type sets of Corona Product of Graphs}
	\begin{theorem}\label{bouncorgraph}
		Let $G$ and $H$ be two simple graphs. Then
		\begin{enumerate}
			\item $\displaystyle{Ecc(G\odot H)=\bigcup_{u_i\in Ecc(G)} V(H_i)}$.
			\item $Ct(G\odot H)=\bigcup\limits_{i=1}^n V(H_i)$.
			\item $\partial(G\odot H)=\bigcup\limits_{i=1}^n V(H_i)$.
		\end{enumerate}
	\end{theorem}
	\begin{proof}
	Let $G$ and $H$ be two simple graphs with vertex sets $\{u_1,u_2,  \ldots ,u_n\}$ and $\{v_1,v_2,  \ldots ,v_m\}$ respectively.\\
	\begin{enumerate}
	   	\item Let $\displaystyle{M=\bigcup_{u_i\in Ecc(G)} V(H_i)}$ and let $w\in Ecc(G\odot H)$. Then either $w=u_i \text{ or }w=v_r^i$.	Suppose $w=u_i$ then there exists $y\in V(G\odot H)$ such that $ecc_{G\odot H}(y)=d_{G\odot H}(y,w)$.\\
		\textbf{Case 1:} $y=u_j$.\\
		Then $\ecc_{G\odot H}(u_j)=d_{G}(u_j,u_i)$. Consider an arbitrary vertex $v_t^i\in V(H_i)$. Then $d_{G\odot H}(u_j,v_t^i)=d_{G}(u_j,u_i)+1$. This implies that $d_{G\odot H}(u_j,u_i)<d_{G\odot H}(u_j,v_t^i)$, a contradiction since $u_i$ is an eccentric vertex of $u_j$.\\
		\textbf{Case 2:} $y=v_t^j$.\\
		Then $\ecc_{G\odot H}(v_t^j)=d_{G\odot H}(v_t^j,u_i)$. Thus $\ecc_{G\odot H}(v_t^j)=d_{G}(u_j,u_i)+1$. Consider an arbitrary vertex $v_r^i\in V(H_i)$. Then $d_{G\odot H}(v_t^j,v_r^i) = d_G(u_j,u_i)+2$. This implies that $d_{G\odot H}(v_t^j,u_i) < d_{G\odot H}(v_t^j,v_r^i)$, a contradiction since $u_i$ is an eccentric vertex of $v_t^j$. Thus $w \neq u_i$.\\
  Hence $w=v_r^i$ and  $v_r^i \in Ecc(G\odot H)$. So, there exist $y\in V(G\odot H)$ such that $ecc_{G\odot H}(y)=d_{G\odot H}(y,v_r^i)$.\\		
		\textbf{Case 1:} $y=u_j$.\\
		$\ecc_{G\odot H}(u_j)=\ecc_G(u_j)+1$ and $d_{G\odot H}(u_j,v_r^i)=d_G(u_j,u_i)+1$. Since $\ecc_{G\odot H}(y)=d_{G\odot H}(y,v_r^i)$, $\ecc_{G\odot H}(u_j)=d_{G\odot H}(u_j,v_r^i)$. This gives $\ecc_G(u_j)=d_G(u_j,u_i)$. Hence $u_i \in Ecc(G)$. Thus $\displaystyle{v_r^i \in M=\bigcup_{u_i\in Ecc(G)} V(H_i)}$.\\
		\textbf{Case 2:} $y=v_s^j$.\\
		$\ecc_{G\odot H}(v_s^j)=\ecc_G(u_j)+2$ and $d_{G\odot H}(v_s^j,v_r^i)=d_G(u_j,u_i)+2$. Since $\ecc_{G\odot H}(y)=d_{G\odot H}(y,v_r^i)$, $\ecc_{G\odot H}(v_s^j)=d_{G\odot H}(v_s^j,v_r^i)$. Thus $\ecc_G(u_j)=d_G(u_j,u_i)$. Hence $u_i \in Ecc(G)$. Thus $\displaystyle{v_r^i \in M=\bigcup_{u_i\in Ecc(G)} V(H_i)}$.	Therefore  $\displaystyle{Ecc(G\odot H)=\bigcup_{u_i\in Ecc(G)} V(H_i)}$. \\
  Next suppose $w\in M$. Then $w=v_r^i$ such that $u_i\in Ecc(G)$.	Since $u_i\in Ecc(G)$ there exist some vertex $u_j \in V(G)$ such that $\ecc_G(u_j)=d_G(u_j,u_i)$. Also $d_{G\odot H}(u_j,v_r^i)=d_G(u_j,u_i)+1$. Then $d_{G\odot H}(u_j,v_r^i)=\ecc_{G\odot H}(u_j)$.	Thus $v_r^i \in Ecc(G\odot H)$. This implies that $M \subseteq Ecc(G\odot H)$. Therefore $\displaystyle{Ecc(G\odot H)=\bigcup_{u_i\in Ecc(G)} V(H_i)}$.\\
		
	\item Let $w\in Ct(G\odot H)$, then either $w= u_i \text{ or }w=v_r^i$. Suppose $w=u_i$. Since $u_i \in Ct(G\odot H)$, $\ecc_{G\odot H}(y)\le \ecc_{G\odot H}(u_i) \text{ for every } y\in N(u_i)$. Consider the vertex $v_r^i\in V(H_i)$, also $v_r^i \in N(u_i)$. $\ecc_{G\odot H}(v_r^i)=ecc_G(u_i)+2$ and $\ecc_{G\odot H}(u_i)=\ecc_G(u_i)+1$. Thus $\ecc_{G\odot H}(v_r^i)>\ecc_{G\odot H}(u_i), v_r^i \in N(u_i)$, a contradiction since $u_i\in Ct(G\odot H)$. Hence $w=v_r^i$ and this implies $w\in V(H_i) \text{ for some }i$. Therefore $Ct(G\odot H)\subseteq \bigcup\limits_{i=1}^nV(H_i)$. Conversely suppose that $w\in \bigcup\limits_{i=1}^nV(H_i)$ then  $w=v_r^i \in V(H_i)$ for some $i$. Neighbouring vertices of $v_r^i \text{ are } N(v_r^i)=\{u_i\}\bigcup \{v_t^i: v_t \in N(v_r)\}$. Now $\ecc_{G\odot H}(v_r^i)=\ecc_G(u_i)+2$. Also $\ecc_{G\odot H}(u_i)=\ecc_G(u_i)+1$ and $\ecc_{G\odot H}(v_t^i)=\ecc_G(u_i)+2$.	Thus  $\ecc_{G\odot H}(y)\leq \ecc_{G\odot H}(v_r^i) \text{ for every }y\in N(v_r^i)$.
		Hence $v_r^i \in Ct(G\odot H)$.	Therefore $Ct(G\odot H)=\bigcup \limits_{i=1}^nV(H_i)$.\\
		
		\item Let $w\in \partial(G\odot H)$, then either $w=u_i \text{ or }w=v_r^i$. Suppose $w=u_i$. Since $u_i\in \partial(G\odot H)$, there exist some vertex $x \in V(G\odot H)$ such that $d_{G\odot H}(x, y)\leq d_{G\odot H}(x, u_i) \text{ for every }y\in N(u_i)$. Consider the vertex $v_r^i \in V(G\odot H)$. Also $v_r^i\in N(u_i)$. Then $d_{G\odot H}(x,v_r^i) \leq d_{G\odot H}(x,u_i)$.\\
		\textbf{Case 1:} $x=u_j$.\\
		Now $d_{G\odot H}(u_j,v_r^i)= d_G(u_i,u_j)+1$. Then $d_{G\odot H}(u_j,v_r^i)> d_{G\odot H}(u_j,u_i)$ which is a contradiction since $u_i \in \partial(G\odot H)$.\\
		\textbf{Case 2:} $x=v_t^j$.\\
		Since $d_{G\odot H}(v_t^j,v_r^i)=d_G(u_j,u_i)+2$ and  $d_{G\odot H}(v_t^j,u_i)=d_G(u_j,u_i)+1$, $d_{G\odot H}(v_t^j,v_r^i)>d_{G\odot H}(v_t^j,u_i)$. This gives a contradiction since $u_i \in \partial(G\odot H)$. Thus $w =v_r^i$, $w\in V(H_i)$. Hence $\partial(G\odot H) \subseteq \bigcup \limits_{i=1}^nV(H_i)$.	Conversely suppose that $v_r^i \in \bigcup \limits_{i=1}^nV(H_i)$. Then $v_r^i \in V(H_i)\text{ for some }i$.\\
		Neighbouring vertices of $v_r^i$ are $N(v_r^i)=\{u_i\}\bigcup \{v_t^i:v_t\in N(v_r)\}$. Now $d_{G\odot H}(u_i,v_r^i)=1$ and $d_{G\odot H}(u_i,v_t^i)=1$. Then there exist vertex $u_i \in V(G\odot H)$ such that $d_{G\odot H}(u_i,y)\leq d_{G\odot H}(u_i,v_r^i)$ for every $y \in N(v_r^i)$. Thus $v_r^i \in \partial (G\odot H)$. Therefore $\partial(G\odot H)=\bigcup \limits_{i=1}^n V(H_i)$.
		\end{enumerate}
	\end{proof}
	\section{Corona Product of Directed Graphs}
	The following is an analogue to the definition of corona product of graphs in \cite{iswadi2011metric}.  Let $D$ be a strong digraph of order $n$ and $H$ be any digraph with $|V(H)|\geq 2$. Then the corona product $D\odot H$ of $D$ and $H$ is defined as a digraph which is formed by taking $n$ copies $H_1,H_2, \ldots, H_n$ of $H$ and connecting  $i^\text{th}$
		vertex of $D$ to the vertices of $H_i$ with a bidirectional edge. Throughout this section $H_i$ is referred to as  $i^\text{th}$ copy of $H$ connected to  $i^\text{th}$ vertex of $D$ in $D\odot H$ for every $i\in \{1,2, \ldots, n\}$.
	Let $V(D)=\{u_1,u_2, \ldots, u_n\}$ and $V(H)=\{v_1,v_2, \ldots, v_m\}$ and $H_i, i\in I=\{1, 2, \ldots, n\}$ be the $i^\text{th}$ copy of $H$. Each vertex $v_r$ of a copy of $H$ attached to a vertex $u_i$ can be uniquely described as $v_r^i$. Thus $V(H_i)=\{v_1^i,v_2^i,  \ldots, v_m^i\}$. Corona product of $D$ and $H$ is a digraph $D\odot H$ such that $V(D\odot H)=V(D)\bigcup \bigg\{\bigcup\limits_{i=1}^n V(H_i)\bigg\}$ and $E(D\odot H)=E(D)\bigcup \bigg\{\bigcup\limits_{i=1}^n E(H_i)\bigg\} \bigcup \bigg\{\bigcup \limits_{i=1}^n\big \{(u_i,v_r^i), (v_r^i,u_i)\big \}, r=1,2, \ldots,m\bigg\}$.	
	\subsection{Distance between two vertices in Corona Product of Digraphs}
	The following proposition gives the distance between two vertices in Corona product.
	\begin{proposition} \label{distcordigraph}
		Let $D$ and $H$ be two strong digraphs with vertex sets $\{u_1,u_2,  \ldots ,u_n\}$ and $\{v_1,v_2, \ldots, v_m\}$ respectively. Then
		\begin{equation*}
			d_{D \odot H}(x,y)  = \begin{cases}
				d_D(u_i,u_j) &\quad x=u_i, y= u_j\\
				d_D(u_i,u_j)+2  &\quad x=u_i, y=v_s^j\\
				d_D(u_i,u_j)+4 &\quad x=v_r^i, y=v_s^j, i\neq j\\
				\min\{d_H(v_r,v_s),4\} &\quad x=v_r^i, y=v_s^j, i=j .
			\end{cases}	
		\end{equation*}
	\end{proposition}
	\begin{proof}
		Given $D$ and $H$ be two  strong digraphs with vertex sets $\{u_1,u_2, \ldots,u_n\}$ and $\{v_1,v_2, \ldots,v_m\}$. Consider the vertices $x,y \in V(D\odot H)$.\\
		\textbf{Case 1:} $x=u_i, y=u_j$. $x, y \in V(D)$.\\
		Then $d_{D\odot H}( x,y)=d_{D}(u_i,u_j)$. Since $u_i\text{ and }u_j$ are vertices of $D$, $d_{D\odot H}(u_i,u_j)=d_D(u_i,u_j)$.\\
		\textbf{Case 2:} $x=u_i, y=v_s^j$. $x\in V(D), y\in V(H_j)$.\\
		Then $d_{D\odot H}(u_i,v_s^j)=\overrightarrow{d}_{D\odot H}(u_i,v_s^j)+	\overrightarrow{d}_{D\odot H}(v_s^j,u_i)$. $d_{D\odot H}(u_i,v_s^j)=\overrightarrow{d}_{D\odot H}(u_i,u_j)+\overrightarrow{d}_{D\odot H}(u_j,v_s^j)+\overrightarrow{d}_{D\odot H}(v_s^j,u_j)+\overrightarrow{d}_{D\odot H}(u_j,u_i)$. Hence $d_{D\odot H}(u_i,v_s^j)=d_D(u_i,u_j)+2$.\\
		\textbf{Case 3:} $x=v_r^i, y=v_s^j, i\neq j$. $x\in V(H_i),y\in V(H_j)$.\\
		$d_{D\odot H}(v_r^i, v_s^j)=\overrightarrow{d}_{D\odot H}( v_r^i,v_s^j)+\overrightarrow{d}_{D\odot H}(v_s^j, v_r^i)$. Then $d_{D\odot H}(v_r^i, v_s^j)=\overrightarrow{d}_{D\odot H}(v_r^i,u_i)+\overrightarrow{d}_{D\odot H}(u_i,u_j)+\overrightarrow{d}_{D\odot H}(u_j,v_s^j)+\overrightarrow{d}_{D\odot H}(v_s^j,u_j)+\overrightarrow{d}_{D\odot H}(u_j,u_i)+\overrightarrow{d}_{D\odot H}(u_i,v_r^i)$. Hence $d_{D\odot H}(v_r^i, v_s^j)=d_D(u_i,u_j)+4$.\\
		\textbf{Case 4:} $x=v_r^i, y=v_s^j, i=j $.\\
		Here the vertices $v_r^i$ and $v_s^j$ lie in the same component. If there exists a bidirectional edge joining these two vertices then ${d}_{D\odot H}(v_r^i,v_s^j)=2$ and if they lie in a directed cycle of three vertices then  ${d}_{D\odot H}(v_r^i,v_s^j)=3$. Otherwise the path $v_r^iu_iv_s^j$ is a path of length $4$ and the length of the shortest $v_r^i-v_s^j$ path is 4. Therefore $d_{D\odot H}(v_r^i, v_s^j)=\min\{d_H(v_r,v_s),4\}$ when $i=j$.
	\end{proof}
	
	\begin{lemma}
		Let $D$ and $H$ be two strong digraphs. Then
		\begin{equation*}
			\ecc_{D\odot H}(x)= \begin{cases}
				\ecc_D(u_i)+2, \quad x=u_i\\
				\ecc_D(u_i)+4, \quad x=v_r^i .
			\end{cases}
		\end{equation*}
	\end{lemma}
	\begin{proof}
		Given $D$ and $H$ be two strong digraphs with vertex sets $\{u_1, u_2,  \ldots ,u_n \}$ and $\{v_1,v_2,  \ldots ,v_m\}$.\\Let $\ecc_{D\odot H}(x)=d_{D\odot H}(x, y) \text{ where } x, y\in V(D\odot H)$.\\
		\textbf{Case 1:} $x= u_i$.\\
		Then $\ecc_{D\odot H}(u_i)=d_{D\odot H}(u_i, y)$. Then either $y \in V(D) \text{ or } y\in V(H_j)$ for some $j \in \{1,2,  \ldots, n\}$.	Suppose $y\in V(D)$ and  $y=u_j$. Then $\ecc_{D\odot H}(u_i)=d_{D\odot H}(u_i,u_j)$. This implies $\ecc_{D\odot H}(u_i)=d_D(u_i,u_j)$.	Consider an arbitrary vertex $v_r^i \in V(H_i)$. Then $d_{D\odot H}(u_i,v_r^i)=d_D(u_i,u_j)+2$. This implies $ecc_{D\odot H}(u_i) < d_{D\odot H}(u_i,v_r^i)$, a contradiction. Thus $y\in V(H_j)$ for some $j \in \{1,2,  \ldots, n\}$. Let $y= v_s^j$. Then $\ecc_{D\odot H}(u_i)=d_{D\odot H}(u_i,v_s^j)$. Thus $\ecc_{D\odot H}(u_i)=d_D(u_i,u_j)+2$. Suppose $\ecc_D(u_i)\neq d_D(u_i,u_j)$, then there exists some vertex $u_t \in V(D)$ such that $d_D(u_i,u_j)< d_D(u_i,u_t)$. Let $v_s^t \in V(H_t)$ be an arbitrary vertex. Since $d_{D\odot H}(u_i,v_s^t))=d_D(u_i,u _t)+2$ and $d_{D\odot H}(u_i,v_s^j)=d_D(u_i,u_j)+2$, $d_{D\odot H}(u_i,v_s^j)<d_D(u_i,v_s^t)$. This gives a	contradiction since $\ecc_{D\odot H}(u_i)=d_{D\odot H}(u_i,v_s^j)$. Thus $\ecc_D(u_i)=d_D(u_i,u_j)$. Hence $\ecc_{D\odot H}(u_i)=\ecc_D(u_i)+2$.	Therefore $\ecc_{D\odot H}(u_i)=\ecc_D(u_i)+2.$\\
		\textbf{Case 2:} $x=v_r^i$.\\
		Then $\ecc_{D\odot H}(v_r^i)=d_{D\odot H}(v_r^i,y)$.
		Then either $y \in V(D) \text{ or } y\in V(H_j)$ for some $j\in \{1,2,  \ldots, n\}$. Suppose $y\in V(D)$ and  $y=u_j$. Then $\ecc_{D\odot H}(v_r^i)=d_{D\odot H}(v_r^i,u_j)$. Hence $\ecc_{D\odot H}(v_r^i)=d_D(u_i,u_j)+2$. Consider an arbitrary vertex $v_m^j \in V(H_j)$. $	d_{D\odot H}(v_r^i,v_m^j)=d_D(u_i,u_j)+4$. This implies that $ecc_{D\odot H}(v_r^i) < d_{D\odot H}(v_r^i,v_m^j)$, a contradiction. Thus $y\in V(H_j)$. Let $y= v_s^j$. Then $ecc_{D\odot H}(v_r^i)=d_{D\odot H}(v_r^i,v_s^j)$. Hence $ecc_{D\odot H}(v_r^i)=d_D(u_i,u_j)+4$. Suppose $\ecc_D(u_i)\neq d_D(u_i,u_j)$. Similarly as in Case 1 this will contradict the fact that $\ecc_{D\odot H}(v_r^i)=d_{D\odot H}(v_r^i,v_s^j)$. Thus $\ecc_D(u_i)=d_D(u_i,u_j)$. Therefore $\ecc_{D\odot H}(v_r^i)=\ecc_D(u_i,u_j)+4$.
	\end{proof}

	\subsection{Radius, Diameter and Center of Corona Product of Digraphs}
	\begin{lemma}
		Let $D$ and $H$ be two strong digraphs. Then $\rad(D\odot H)=\rad(D)+2$.
	\end{lemma}
	\begin{proof}
		Given $D$ and $H$ be two strong digraphs. Consider the vertices $v_r^i, u_i \in V(D\odot H)$. Then $\ecc_{D\odot H}(v_r^i)=\ecc_D(u_i)+4, v_r^i\in V(H_i)$ and $\ecc_{D\odot H}(u_i)=\ecc_D(u_i)+2, u_i\in V(D)$. Thus
		\begin{align*}
			\rad(D\odot H)&=\min\{\ecc_{D\odot H}(x):x\in V(D\odot H)\}\\
			&=\min\{\ecc_{D\odot H}(u_i)+2:u_i\in V(D)\}\\
			&=\min\{\ecc_{D\odot H}(u_i):u_i\in V(D)\}+2\\
			&=\rad(D)+2.
		\end{align*}
	\end{proof}
	\begin{lemma}
		Let $D$ and $H$ be two strong digraphs. Then $\diam(D\odot H)=\diam(D)+4$.
	\end{lemma}
	\begin{proof}
		Given $D$ and $H$ be two strong digraphs. Consider the vertices $v_r^i, u_i \in V(D\odot H)$. Then $\ecc_{D\odot H}(v_r^i)=\ecc_D(u_i)+4, v_r^i\in V(H_i)$ and $\ecc_{D\odot H}(u_i)=\ecc_D(u_i)+2, u_i\in V(D)$. Thus
		\begin{align*}
			\diam(D\odot H)&=\max\{\ecc_{D\odot H}(x):x\in V(D\odot H)\}\\
			&=\max\{\ecc_{D\odot H}(u_i)+4:u_i\in V(D)\}\\
			&=\max\{\ecc_{D\odot H}(u_i):u_i\in V(D)\}+4\\
			&=\diam(D)+4.
		\end{align*}
	\end{proof}
	\begin{theorem}
		Let $D$ and $H$ be two strong digraphs. Then $C(D\odot H) = C(D)$.
	\end{theorem}
	\begin{proof}
		Let $D$ and $H$ be two strong digraphs with vertex sets $\{u_1,u_2,  \ldots ,u_n\}$ and $\{v_1,v_2,  \ldots ,v_m\}$ respectively. Let $x \in C(D\odot H)$ then  $\ecc_{D \odot H}(x)=\underset{y\in V(D\odot H)}\min \ecc(y)$. Since $x \in V(D\odot H), x \in V(D) \text{ or } x\in V(H_i) \text{ for some } i\in \{1,2,  \ldots, n\}$. Suppose $x=v_r^i \in V(H_i)$. $\ecc_{D\odot H}(v_r^i)=\ecc_D(u_i)+4$ and $\ecc_{D\odot H}(u_i)=\ecc_D(u_i)+2$. This implies that $\ecc_{D\odot H}(u_i) < \ecc_{D\odot H}(v_r^i)$ which is a contradiction since $v_r^i\in C(D\odot H)$. 	Therefore $x \in V(D)$. Thus $x=u_i \text{ for some }i \in \{1,2,  \ldots ,n\}$. Suppose $\ecc_D(u_i)\neq \rad(D)$ then there exists some vertex $u_j \in V(D)$ such that $\ecc_D(u_j)<\ecc_D(u_i)$. Since $\ecc_{D\odot H}(u_j)=\ecc_D(u_j)+2$ and $\ecc_{D\odot H}(u_i)=\ecc_D(u_i)+2$, $\ecc_{D\odot H}(u_j) < \ecc_{D\odot H}(u_i)$. This contradicts the fact that $u_i\in C(D\odot H)$. Thus $\ecc_D(u_i)=\rad(D)$. Hence $u_i \in C(D)$. Now it remains to prove that $C(D) \subseteq C(D\odot H)$.
		Conversely suppose $x \in C(D)$. Let $x=u_i \in C(D)$. Suppose $u_i \notin C(D\odot H)$ then there exists some vertex $y \in V(D\odot H)$ such that $\ecc_{D\odot H}(y)<\ecc_{D\odot H}(u_i)$. Then either $y \in V(D)$ or $y\in V(H_j)$ for some $j \in \{1,2,  \ldots ,n\}$.\\
		\textbf{Case 1 :} $y \in V(H_j)$ for some $j \in \{1,2,  \ldots ,n\}$.\\
		Let $y=v_r^j$. Then $\ecc_{D\odot H}( v_r^j)=\ecc_D(u_j)+4$ and $\ecc_{D\odot H}(u_i)=\ecc_D(u_i)+2$. Since $\ecc_{D\odot H}(y)<\ecc_{D\odot H}(u_i)$, $\ecc_D(u_j)<\ecc_D(u_i)$ a contradiction since $u_i \in C(D)$.\\
		\textbf{Case 2 :} $y=u_j \text{ for some }j \in \{1,2,  \ldots,n\}$.\\
		$\ecc_{D\odot H}(u_j)=\ecc_D(u_j)+2$ and $\ecc_{D\odot H}(u_i)=\ecc_D(u_i)+2$. This implies that $\ecc_D(u_j)<\ecc_D(u_i)$ contradiction since $u_i \in C(D)$. Thus $u_i\in C(D\odot H)$. Therefore $C(D\odot H) = C(D)$.
	\end{proof}
	\subsection{Boundary-type sets of Corona Product of Digraphs}
	\begin{theorem}
		\label{boundarycoronadig}
		Let $D$ and $H$ be two strong digraphs. Then 
		\begin{enumerate}
			\item $\displaystyle{Per(D\odot H)=\bigcup_{u_i\in Per(D)} V(H_i)}$.
			\item $\displaystyle{Ecc(D\odot H)=\bigcup_{u_i\in Ecc(D)} V(H_i)}$.
			\item $Ct(D\odot H)=\bigcup\limits_{i=1}^n V(H_i)$.
			\item $\partial(D\odot H)=\bigcup\limits_{i=1}^n V(H_i)$.	
		\end{enumerate}
	\end{theorem}
	
	\begin{proof}
		Let $D$ and $H$ be two strong digraphs with vertex sets $\{u_1,u_2,  \ldots ,u_n\}$ and $\{v_1,v_2,  \ldots ,v_m\}$ respectively.\\
		\begin{enumerate}
			\item Let $\displaystyle{M=\bigcup_{u_i\in Per(D)} V(H_i)}$.
		Suppose $w \in Per(D\odot H)$ then $\ecc_{D\odot H}(w)=\diam(D\odot H)$. Since $w \in Per(D\odot H)$, $w \in V(D\odot H)$ this implies $w\in V(D) \text{ or }w\in V(H_i)$ for some $i\in \{1,2,  \ldots ,n\}$. Suppose $w= u_i \in V(D)$ and let $v_s^i$ be arbitrary vertex in $V(H_i)$. Then $ v_s^i\in V(D\odot H)$. Now
		$\ecc_{D\odot H}( v_s^i)=\ecc_D(u_i)+4$ and $\ecc_{D\odot H}(u_i)=\ecc_D(u_i)+2$. Since $\ecc_{D\odot H}(w)=\ecc_{D\odot H}(u_i)$, $\ecc_{D\odot H}(w)<\ecc_{D\odot H}( v_s^i)$ a contradiction since $w \in Per(D\odot H)$. Thus $w\in V(H_i)$ for some $i \in\{1,2,\ldots,n\}$. Let $w=v_r^i \in V(H_i)$.	Suppose  $u_i\notin Per(D)$. Then there exists some vertex $u_m\in V(D)$ such that $\ecc_D(u_i)<\ecc_D(u_m)$. Since $u_m \in V(D), v_r^m \in V(D\odot H)$. Also 	$\ecc_{D\odot H}( v_r^m)=\ecc_D(u_m)+4$ and $\ecc_{D\odot H}( v_r^i)=\ecc_D(u_i)+4$. Since $\ecc_D(u_i)<\ecc_D(u_m)$, $\ecc_{D\odot H}( v_r^i)<\ecc_{D\odot H}( v_r^m)$. This is a contradiction since $v_r^i\in Per(D\odot H)$. Thus $u_i\in Per(D)$. Hence $u_i \in Per(D)$.\\ Conversely let $\displaystyle{w \in M=\bigcup_{u_i\in Per(D)} V(H_i)}$. Then $w=v_r^i \text{ and } u_i\in Per(D)$. Then $\ecc_D(u_i)=\diam(D)$. Suppose $\ecc_{D\odot H}( v_r^i)\neq diam(D\odot H)$ then there exists some vertex $v_s^j\in V(D\odot H) \text{ such that } \ecc_{D\odot H}( v_r^i)<\ecc_{D\odot H}(v_s^j)$. Since $\ecc_{D\odot H}( v_r^i)=\ecc_D(u_i)+4$ and $\ecc_{D\odot H}(v_s^j)=\ecc_D(u_j)+4$, $\ecc_D(u_i)<\ecc_D(u_j)$	contradiction since $\ecc_D(u_i)=\diam(D)$. Thus $\ecc_{D\odot H}( v_r^i)=\diam(D\odot H)$. Hence $ v_r^i\in Per(D\odot H)$.  Therefore $\displaystyle{Per(D\odot H)=\bigcup_{u_i\in Per(D)} V(H_i)}$.\\
		
		\item Let $\displaystyle{M=\bigcup_{u_i\in Ecc(D)} V(H_i)}$ and let $w\in Ecc( D\odot H)$. Then either $w=u_i \text{ or }w=v_r^i$.  Suppose first $w=u_i$ then there exists $y\in V( D\odot H)$ such that $\ecc_{ D\odot H}(y)=d_{ D\odot H}(y,w)$.\\
				\textbf{Case 1:} $y=u_j$.\\
				Then $\ecc_{ D\odot H}(u_j)=d_{D\odot H}(u_j,u_i)$. Consider an arbitrary vertex $v_t^i\in V(H_i)$. Then $d_{ D\odot H}(u_j,v_t^i)=d_{D}(u_j,u_i) +2$. This implies that $d_{ D\odot H}(u_j,u_i)<d_{ D\odot H}(u_j,v_t^i)$, a contradiction to fact that $u_i$ is an eccentric vertex of $u_j$.\\
				\textbf{Case 2:} $y=v_t^j$.\\
				Then $\ecc_{ D\odot H}(v_t^j)=d_{ D\odot H}(v_t^j,u_i)$. Thus $\ecc_{ D\odot H}(v_t^j)=d_{ D}(u_j,u_i) +2$. Consider an arbitrary vertex $v_r^i\in V(H_i)$. Then $d_{ D\odot H}(v_t^j,v_r^i) = d_ D(u_j,u_i) +4$. This implies that $d_{ D\odot H}(v_t^j,u_i) < d_{ D\odot H}(v_t^j,v_r^i)$, a contradiction since $u_i$ is an eccentric vertex of $v_t^j$. From both cases $w \neq u_i$. Hence $w=v_r^i$. So $v_r^i \in Ecc( D\odot H)$. Then there exist $y\in V( D\odot H)$ such that $	\ecc_{ D\odot H}(y)=d_{ D\odot H}(y,v_r^i)$.\\		
				\textbf{Case 3:} $y=u_j$.\\
				$\ecc_{ D\odot H}(u_j)=\ecc_ D(u_j) +2$ and $d_{ D\odot H}(u_j,v_r^i)=d_ D(u_j,u_i) +2$. Since $\ecc_{ D\odot H}(y)=d_{ D\odot H}(y,v_r^i)$, $\ecc_{ D\odot H}(u_j)=d_{ D\odot H}(u_j,v_r^i)$. This gives $\ecc_ D(u_j)=d_ D(u_j,u_i)$. Hence $u_i \in Ecc( D)$. Thus $v_r^i \in M=\{v_r^i:u_i\in Ecc( D)\}$.\\
				\textbf{Case 4:} $y=v_s^j$.\\
				$\ecc_{ D\odot H}(v_s^j)=\ecc_ D(u_j) +4$ and $d_{ D\odot H}(v_s^j,v_r^i)=d_ D(u_j,u_i) +4$. Since $\ecc_{ D\odot H}(y)=d_{ D\odot H}(y,v_r^i)$, $\ecc_{ D\odot H}(v_s^j)=d_{ D\odot H}(v_s^j,v_r^i)$. Thus $\ecc_ D(u_j)=d_ D(u_j,u_i)$. Hence $u_i \in Ecc( D)$. Thus $v_r^i \in M=\{\cup V(H_i):u_i\in Ecc( D)\}$.	Therefore  $Ecc( D\odot H)\subseteq M$. Suppose $w\in M$. Then $w=v_r^i$ such that $u_i\in Ecc( D)$.	Since $u_i\in Ecc( D)$ there exist some vertex $u_j \in V( D)$ such that $\ecc_ D(u_j)=d_ D(u_j,u_i)$. Also $d_{ D\odot H}(u_j,v_r^i)=d_ D(u_j,u_i) +2$. Then $d_{ D\odot H}(u_j,v_r^i)=\ecc_{ D\odot H}(u_j)$.	Thus $v_r^i \in Ecc( D\odot H)$. This implies that $M \subseteq Ecc( D\odot H)$. Therefore $\displaystyle{Ecc(D\odot H)=\bigcup_{u_i\in Ecc(D)} V(H_i)}$.\\
		
		\item Let $w\in Ct( D\odot H)$, then either $w= u_i \text{ or }w=v_r^i$. Suppose $w=u_i$. Since $u_i \in Ct( D\odot H)$, $\ecc_{ D\odot H}(y)\le \ecc_{ D\odot H}(u_i) \text{ for every } y\in N(u_i)$. Consider the vertex $v_r^i\in V(H_i)$, also $v_r^i \in N(u_i)$. $\ecc_{ D\odot H}(v_r^i)=\ecc_ D(u_i) +4$ and $\ecc_{ D\odot H}(u_i)=\ecc_ D(u_i) +2$. Thus $\ecc_{ D\odot H}(v_r^i)>\ecc_{ D\odot H}(u_i), v_r^i \in N(u_i)$, a contradiction since $u_i\in Ct( D\odot H)$. Hence $w=v_r^i$ and this implies $w\in V(H_i) \text{ for some }i$. Therefore $Ct( D\odot H)\subseteq \bigcup\limits_{i=1}^nV(H_i)$. Conversely suppose that $w\in \bigcup\limits_{i=1}^nV(H_i)$ then  $w=v_r^i \in V(H_i)$ for some $i$. Neighbouring vertices of $v_r^i \text{ are } N(v_r^i)=\{u_i\}\bigcup \{v_t^i: v_t \in N(v_r)\}$. Now $\ecc_{ D\odot H}(v_r^i)=\ecc_ D(u_i) +4$. Also $\ecc_{ D\odot H}(u_i)=\ecc_ D(u_i) +2$ and $\ecc_{ D\odot H}(v_t^i)=\ecc_ D(u_i) +4$.	Thus  $\ecc_{ D\odot H}(y)\leq \ecc_{ D\odot H}(v_r^i) \text{ for every }y\in N(v_r^i).$
				Hence $v_r^i \in Ct( D\odot H)$.	Therefore $Ct( D\odot H)=\bigcup \limits_{i=1}^nV(H_i)$.\\ 
				
			\item  Let $w\in \partial( D\odot H)$, then either $w=u_i \text{ or }w=v_r^i$. Suppose $w=u_i$. Since $u_i\in \partial( D\odot H)$, there exist some vertex $x \in V( D\odot H)$ such that $d_{ D\odot H}(x, y)\leq d_{ D\odot H}(x, u_i) \text{ for every }y\in N(u_i)$. Consider the vertex $v_r^i \in V( D\odot H)$. Also $v_r^i\in N(u_i)$. Then $d_{ D\odot H}(x,v_r^i) \leq d_{ D\odot H}(x,u_i)$.\\
				\textbf{Case 1:} $x=u_j$.\\
				Now $d_{ D\odot H}(u_j,v_r^i)= d_ D(u_i,u_j) +2$. Then $d_{ D\odot H}(u_j,v_r^i)> d_{ D\odot H}(u_j,u_i)$ which is a contradiction since $u_i \in \partial( D\odot H)$.\\
				\textbf{Case 2:} $x=v_t^j$.\\
				Since $d_{ D\odot H}(v_t^j,v_r^i)=d_ D(u_j,u_i) +4$ and  $d_{ D\odot H}(v_t^j,u_i)=d_ D(u_j,u_i) +2$, $d_{ D\odot H}(v_t^j,v_r^i)>d_{ D\odot H}(v_t^j,u_i)$. This gives a contradiction since $u_i \in \partial(D\odot H)$. Thus $w =v_r^i$, $w\in V(H_i)$. Hence $\partial(D\odot H) \subseteq \bigcup \limits_{i=1}^nV(H_i)$.	Conversely, suppose that $v_r^i \in \bigcup \limits_{i=1}^nV(H_i)$. Then $v_r^i \in V(H_i)\text{ for some }i$.\\
				Neighbouring vertices of $v_r^i$ are $N(v_r^i)=\{u_i\}\bigcup \{v_t^i:v_t\in N(v_r)\}$. Now $d_{D\odot H}(u_i,v_r^i)=2$ and $d_{G\odot H}(u_i,v_t^i)=2$. Then there exist vertex $u_i \in V(D\odot H)$ such that $d_{G\odot H}(u_i,y)\leq d_{D\odot H}(u_i,v_r^i)$ for every $y \in N(v_r^i)$. Thus $v_r^i \in \partial (G\odot H)$. Therefore $\partial(D\odot H)=\bigcup \limits_{i=1}^n V(H_i)$.	
  \end{enumerate}
	\end{proof}
\section*{Conflict of Interest}
The authors declare no conflict of interest.
	\bibliographystyle{amsplain}
	\bibliography{bibliography}
\end{document}